\title{Dense graphs are antimagic}

\author{
N. Alon
\thanks{
Department of Mathematics,
Tel Aviv University, Tel Aviv 69978, Israel. E--mail:
nogaa@post.tau.ac.il.
Research supported in part
by a USA Israeli BSF grant, by a grant from
the Israel Science Foundation
and by the Hermann Minkowski Minerva
Center for Geometry at Tel Aviv University.
}
\and
G. Kaplan
\thanks{
Department of Computer Sciences
The Academic College of Tel-Aviv-Yaffo
Tel-Aviv 61161, Israel and
Department of Mathematics
School of Mathematical Sciences
Tel Aviv University, Tel Aviv 69978, Israel.
}
\and
A. Lev
\thanks{
Department of Computer Sciences
The Academic College of Tel-Aviv-Yaffo
Tel-Aviv 61161, Israel and
Department of Mathematics
School of Mathematical Sciences
Tel Aviv University, Tel Aviv 69978, Israel.
}
\and
Y. Roditty
\thanks{
School of Computer Sciences,
Tel Aviv University, Tel Aviv 69978, Israel and
Department of Computer Sciences, The Academic College of
Tel-Aviv-Yaffo, Tel-Aviv 61161, Israel.
email: jr@post.tau.ac.il}
\and
R. Yuster
\thanks{Department of Mathematics, University of
Haifa-Oranim, Tivon 36006,
Israel.
email: raphy@research.haifa.ac.il}
}

\date{} 

\documentstyle [11pt]{article}
\newtheorem{theorem}{Theorem}[section]

\newtheorem{lemma}[theorem]{Lemma}

\newtheorem{conj}[theorem]{Conjecture}

\newcommand{\ignore}[1]{}
\def\square{\vrule height6pt width7pt depth1pt}
\def\endpf{\hfill\square\bigskip}

\begin{document}
\maketitle
\setcounter{page}{1}
\begin{abstract}
An {\em antimagic labeling} of a graph with $m$ edges and $n$ vertices is a
bijection from the set of edges to the integers $1,\ldots,m$ such
that all $n$ vertex sums are pairwise distinct, where
a vertex sum is the sum of
labels of all edges incident with the same vertex.
A graph is called {\em antimagic} if it has an antimagic labeling.
A conjecture of Ringel (see \cite{HaRi}) states that every connected graph, but $K_2$, is antimagic.
Our main result validates this conjecture for graphs
having minimum degree $\Omega(\log n)$.
The proof combines probabilistic arguments with simple 
tools from analytic number
theory and combinatorial techniques.
We also prove that complete partite graphs (but $K_2$)
and graphs with maximum degree at least $n-2$ are antimagic.\\
{\bf AMS classification code:} 05C78
\end{abstract}

\section{Introduction}
All graphs in this paper are finite, undirected and simple.
We follow the notation and terminology of \cite{Bo}.
An {\em antimagic labeling} of a graph with $m$ edges and $n$ vertices is a
bijection from the set of edges to the integers $1,\ldots,m$ such
that all $n$ vertex sums are pairwise distinct, where a 
vertex sum is the sum of
labels of all edges incident with the same vertex.
A graph is called {\em antimagic} if it has an antimagic labeling.
The following is conjectured in \cite{HaRi}:
\begin{conj}
\label{c1}
Every connected graph, but $K_2$, is antimagic.
\end{conj}
In this paper we prove that conjecture \ref{c1} holds for several 
classes of graphs.
Our main result validates Conjecture \ref{c1} for all graphs with minimum degree $\Omega(\log n)$.
\begin{theorem}
\label{t1}
There exists an absolute constant $C$ such that every graph with $n$ vertices
and minimum degree at least $C \log n$ is antimagic.
\end{theorem}
In fact, our proof can be optimized to obtain that even
a lower bound of $\Omega(\log n/\log\log n)$ for the
minimum degree suffices. 
Since this improvement in negligible and complicates the computations
significantly, we omit its proof. The proof of Theorem \ref{t1} 
requires several probabilistic tools and some simple techniques
from analytic number theory,
together with some combinatorial ideas, 
and is presented in the next two sections.

It is rather straightforward to prove that $n$-vertex graphs with a vertex of degree $n-1$
are antimagic. It is somewhat surprising that this ceases to be straightforward when the minimum degree is $n-2$.
The following theorem is proved in Section 4.
\begin{theorem}
\label{t2}
If $G$ has $n \geq 4$ vertices and $\Delta(G) \geq n-2$ then $G$ is antimagic.
\end{theorem}
It is still an open problem to decide whether connected 
graphs with $\Delta(G) \geq n-k$ and $n >n_0(k)$ are antimagic,
for any fixed $k \geq 3$.

In the final section we prove the following.
\begin{theorem}
\label{t3}
All complete partite graphs, but $K_2$, are antimagic.
\end{theorem}

\section{A few probabilistic lemmas}
In the next section we describe some lemmas that will be useful in the
proof of Theorem \ref{t1}. 
\begin{lemma}
\label{l21}
There are absolute positive constants $c_1, c_2 >0$ 
such that the following holds.
Let $t$ be a positive integer, let $A \subset \{1,2, \ldots ,t\}$
be a set of size $|A| \geq t-2d$, where $d$ is a positive integer,
$d \leq t/30$, and define $p =\lfloor td^{1/2} \rfloor$.
Let $w=e^{2\pi i/p}$ be a primitive root of
$1$ of order $p$,  and let $a_1,a_2$ be two distinct randomly chosen
elements of $A$ (where all pairs are equally likely).

\noindent
(i) If $x$ is an integer, and 
$0<x < d^{1/2}$ or $p-d^{1/2} < x < p$, then with probability at least
$c_1$, 
\begin{equation}
\label{e1}
|\frac{w^{a_1x} +w^{a_2x}}{2}| \leq 1-c_2 \frac{min(x,p-x)^2}{d}.
\end{equation}

\noindent
(ii) If $x$ is an integer and $d^{1/2} \leq x \leq p-d^{1/2}$ then with
probability at least $c_1$,
\begin{equation}
\label{e2}
|\frac{w^{a_1x} +w^{a_2x}}{2}| \leq 1-c_2.
\end{equation}
\end{lemma}
{\bf Proof}\,  
Note that
$$
|\frac{w^{a_1x} +w^{a_2x}}{2}|=|\frac{1 +w^{(a_2-a_1)x}}{2}|
=|cos (\frac{(a_2-a_1) \pi x}{p})|.
$$
If $0 <x < d^{1/2}$ then, since $A$ contains most of the numbers in
the interval $\{1,2, \ldots ,t\}$, it follows that with probability 
$\Omega(1)$, $(a_2-a_1)$ is between, say, $t/4$ and $3t/4$,
supplying the assertion of (\ref{e1}), since
$|cos(y)|=1-\Theta(y^2)$ for every $0<y<3\pi/4$.
As $|cos(y)|=|cos(\pi-y)|$ the result for $p-d^{1/2} < x <p$
follows as well.

If $ d^{1/2} \leq x \leq p-d^{1/2}$ then $((a_2-a_1)x)(mod~p)$ is 
between, say, $p/4$ and $3p/4$ with probability $\Omega(1)$,
implying the assertion of (\ref{e2}). \endpf

\begin{lemma}
\label{l22}
There is an absolute positive constant $C_1$ such that the following holds.
Let $t,d$ be integers, $t/30 \geq d \geq  \lfloor C_1 \log t \rfloor$, 
and let $p$ and $w$ be as in
Lemma \ref{l21}. Let $a_{i1},a_{i2}$, $1 \leq i \leq d$, be $2d$ pairwise
distinct elements of $\{1,2, \ldots ,t\}$, chosen randomly with all choices of
$d$ pairwise disjoint pairs being equally likely. 
Define
$$
T(x)=\prod_{i=1}^d \frac{w^{a_{i1}x}+ w^{a_{i2}x}}{2}.
$$
Then with probability at least $1-1/t^2$ the following holds:

\noindent
(i) For every integer $x$ satisfying $0 < x < d^{1/2}$ or
$p-d^{1/2}<x<p$, $|T(x)| \leq e^{-min(x,p-x)^2}$.

\noindent
(ii) For every integer $x$ satisfying $ d^{1/2} \leq x \leq p-d^{1/2}$,
$|T(x)| \leq \frac{1}{t^2}$.
\end{lemma}
{\bf Proof}\, 
Fix an integer $x$ between $1$ and $p-1$.
The random choice of the elements $a_{i1},a_{i2}$ can be obviously
done as follows. Start with $A_1=\{1,2, \ldots ,t\}$, choose a random
pair $a_{11},a_{12}$ of distinct elements of $A_1$, and define $A_2=
A_1-\{a_{11},a_{12}\}$. Next choose a random pair $a_{21},a_{22}$
of distinct elements of $A_2$ and omit them from $A_2$ to get
$A_3$, etc. Since in each step the remaining size of the set $A_i$
is bigger than $t-2d$, we can apply Lemma \ref{l21} and conclude that
in each step, with probability at least $c_1$, the absolute value
of the obtained term
$\frac{w^{a_{i1}x}+w^{a_{i2}x}}{2}$ is at most 
$1-c_2 \frac{min(x,p-x)^2}{d}$
in case $0<x < d^{1/2}$ or $p-d^{1/2} < x < p$, and at most $1-c_2$
in case $d^{1/2} \leq x \leq p-d^{1/2}.$ Let us call the $i$-th chosen
pair $a_{i1},a_{i2}$ {\em successful} if this inequality holds.

By the standard estimates for Binomial distributions (see, e.g., \cite{AlSp}, Appendix A)
and by our assumption that $d=\Omega(\log t)$ it follows that if $C_1$ is chosen appropriately,
then with probability at least $1-1/t^2$, for every fixed admissible
$x$ there are sufficiently many successful pairs to insure the assertion of the lemma. \endpf

\begin{lemma}
\label{l23}
There are absolute positive constants $C_1,C_2$ such that the following
holds.
Let $t$ and $d \geq \lfloor C_1 \log t \rfloor$ be as in Lemma \ref{l22}, and let
$a_{i,1},a_{i,2}$, $1 \leq i \leq d$ be $2d$ pairwise
distinct elements of $\{1,2, \ldots ,t\}$, chosen randomly with all choices of
$d$ pairwise disjoint pairs being equally likely. Then, with probability
at least $1-1/t^2$ the following holds.
For each $i$,
$1 \leq i \leq d$, choose $j_i \in \{1,2\}$  randomly, independently
and uniformly, and consider the random sum $Q=\sum_{i=1}^d a_{i,j_i}$.
Then, for every integer $S$, the probability that $Q$ is equal 
$S$ is at most $\frac{C_2}{td^{1/2}}.$
\end{lemma}
{\bf Proof}\, 
Put $p=\lfloor td^{1/2} \rfloor $. It suffices to prove that the 
probability that $Q$ is equal to $S$ modulo $p$ is at most $O(1/p)$. 
It is thus enough to show that this is the case if the conclusions of
Lemma \ref{l22} hold. Note that the probability that $Q$ is equal to
$S$ modulo $p$ is precisely
$$
\frac{1}{p} \sum_{x=0}^{p-1} \prod_{i=1}^d \frac{w^{a_{i1}x}+
w^{a_{i2}x}}{2} w^{-Sx}= \frac{1}{p}\sum_{x=0}^{p-1} T(x) w^{-Sx}.
$$
However, if the conclusions of Lemma \ref{l22} hold, then the term
corresponding to $x=0$ in the above sum is $1/p$, and the absolute value
of the term corresponding to $x$ is  at most 
$e^{-min(x,p-x)^2}/p$ for each $x$ satisfying 
$0 <x < d^{1/2}$ or $p-d^{1/2} < x <p$, and is at most
$\frac{1}{pt^2}$ for each other $x$ in the sum. As $\sum_{x>0} e^{-x^2} =O(1)$
and as $1/t^2<1/p$, the desired result follows.  \endpf

In the proof of Theorem \ref{t1} we need to use the symmetric
Lov\'asz Local Lemma \cite{ErLo}. Here it is,
following the notations in \cite{AlSp}.
Let $A_1,\ldots,A_s$ be events in an arbitrary probability space.
\begin{lemma}[The Local Lemma, symmetric version]
\label{l24}
If there are positive constants $p$ and $r$ such that $p(r+1) < 1/3$,
each $A_i$ is mutually independent of all other events but at most $r$,
and $\Pr[A_i] \leq p$ for all $i$, 
then with positive probability no event $A_i$ holds. \endpf
\end{lemma}

\section{Proof of Theorem \ref{t1}}
Let $C$ be a sufficiently large absolute constant such that
for all sufficiently large positive integer $n$, and for all $t \in [dn/2,dn]$ where $d=\lfloor C \log n \rfloor$
it holds that $t/30 -1 \geq d \geq \lfloor C_1 \log t \rfloor + 1$ where $C_1$ is the constant from
Lemma \ref{l23}. In order to prove Theorem \ref{t1} it suffices to prove that for all sufficiently large $n$,
if $G$ is a graph with $n$ vertices, $m$ edges, and $\delta(G) \geq d$, then $G$ is antimagic.

It is convenient  to split the description of the proof into  five
phases.

Phase 1: As long as there are two adjacent vertices each having degree at least $d+1$, we assign the edge
connecting them the highest yet unused label and delete the edge.
Let $G'$ denote the spanning subgraph of $G$ obtained at the end of this process.
Denote the set of vertices with degree
$d$ in $G'$ by $A$ and the set of vertices with degree at least $d+1$ in $G'$ by $B$.
Notice that it is possible that $B = \emptyset$ and also note that $B$ induces an independent set in $G'$.
Also, each vertex $v$ has a partial sum denoted $r(v)$, which is the
sum of the labels assigned to the edges incident with $v$ that were deleted in this phase.
Let $t \leq m$ denote the number of edges 
of $G'$ and notice that $t \in [dn/2,dn]$.
Hence, our goal is to assign the set $\{1,\ldots,t\}$ of labels to 
the edges of $G'$
such that all vertex sums (including the partial 
sums contributed from the labels assigned already) are distinct.

Phase 2: We partition the edge-set of $G'$ into $t/2$ pairs 
(clearly, we may assume, without loss of generality, 
that $t$ is even) as follows.
Let $d'(v)$ denote the degree of $v$ in $G'$.
For each $v \in B$ we arbitrarily choose a set $F(v)$ of edges incident with $v$, such that $|F(v)|$ is even and
$$
\left| (d'(v)-|F(v)|) - d \right| \leq 1.
$$
Notice that for any two distinct vertices $v$ and $u$ of $B$ we have $F(v) \cap F(u) = \emptyset$.
We partition each $F(v)$ into pairs, arbitrarily.
Put $k= (t - |\cup_{v \in B} F(v)|)/2$.
We partition the remaining edges of $G'$ not in $\cup_{v \in B} F(v)$ into $k$ pairs,
with the property that the two edges in each pair do not share a common endpoint. This can be done by
observing that the complement of the line graph of the remaining edges has a very high minimum degree and therefore
has a perfect matching.
Finally, for each edge $e \in G'$, let $p(e)$ denote the edge paired with $e$.

Phase 3: We randomly partition the set of $t$ labels into $t/2$ pairs of labels. The randomly selected  label pairs
are then arbitrarily assigned to the edge pairs created in Phase 2. For $e \in G'$ let $L(e)$ denote the pair of labels
assigned to the pair $\{e, p(e)\}$. Clearly $L(p(e))=L(e)$.

Phase 4: For each $v \in B$ let $f(v)$ denote the sum of the labels assigned to the edges of $F(v)$. Notice that although
we have yet to specify which edge gets 
which label (there are two choices), $f(v)$ is well defined.
Now, for each $v \in B$, let $H(v)$ denote the set of edges incident with $v$ in $G'$, and not belonging to $F(v)$.
For each $v \in A$, let $H(v)$ denote the set of edges incident with $v$ in $G'$.
Notice that $d+1 \geq |H(v)| \geq d-1$ for all $v \in V$ (in fact, $|H(v)|=d$ for $v \in A$).
For each $v \in V$, consider the set of pairs of labels $\{ L(e) ~: ~ e \in H(v)\}$.
There are $2^{|H(v)|}$ choices to select one label from each pair $L(e)$, and each choice yields a possible
sum. Denote the set of possible sums by $Q(v)$, and notice that $|Q(v)| \leq 2^{|H(v)|}$.
By Lemma \ref{l23}, with probability at least $1-1/t^2$, no
specific value in $Q(v)$ is obtained by more than a fraction of $C_2/td^{1/2}$ of the $2^{|H(v)|}$ choices
for a suitable absolute constant  $C_2$. We therefore {\em fix} an assignment of label pairs to edge pairs having the
property that for {\em all} $v \in V$, no value from $Q(v)$ is obtained by more than a fraction of
$C_2/td^{1/2}$ of the $2^{|H(v)|}$ possible selections.

Phase 5: For each pair $\{e, p(e)\}$ we flip a coin to decide which edge gets which label from $L(e)$.
All $t/2$ decisions are independent. Notice that the final weight of each $v \in B$ is a random variable given by adding to
$r(v)+f(v)$ a number of $Q(v)$ which corresponds to the random coin flip.
Similarly, the final weight of each $v \in A$ is a random variable given by adding to $r(v)$ a number of $Q(v)$.
We claim that with positive probability, no two vertices of $v$ will end up with the same final weight.
For a pair of vertices $u,v$, let $B(u,v)$ denote the event that both $u$ and $v$ end up with the same final weight.
We need to show that with positive probability no $B(u,v)$ holds. By our arguments from Phase 4,
$$
\Pr[B(u,v)] \leq \frac{C_2}{td^{1/2}} \leq \frac{2C_2}{nd^{3/2}}.
$$
We show that $B(u,v)$ is independent of all other events
but at most $O(nd)$. Indeed, let $Z$ denote the set of vertices
that are endpoints of edges in $H(v) \cup H(u)$ and their matched edges. Clearly $|Z| \leq 6(d+1)+2$.
Now any combination of the events $B(x,y)$ where neither $x$ nor $y$ are in $Z$ is independent of
$B(u,v)$. Thus, $B(u,v)$ is independent of all but 
at most $(6(d+1)+2)n$ other events.
Since
$$
\frac{2C_2}{nd^{3/2}} \cdot ((6(d+1)+2)n+1) << \frac{1}{3}
$$
for $n$ sufficiently large, we get, using Lemma \ref{l24}, that 
with positive probability no $B(u,v)$ holds. This completes 
the proof of Theorem \ref{t1}. \endpf

\section{Proof of Theorem \ref{t2}}
The proof of Theorem \ref{t2} is divided into four lemmas.
The first lemma handles the case of maximum degree $n-1$.
\begin{lemma}
\label{l41}
If $G$ has $n$ vertices and $\Delta(G)=n-1$ then $G$ is antimagic.
\end{lemma}
{\bf Proof}\,
Assume $G$ has $m$ edges, and let $v$ be a vertex of degree $n-1$.
Assign the distinct labels $1,\ldots,m-n+1$ arbitrarily to all the $m-n+1$
edges not incident with $v$. Denote the $n-1$ neighbors of $v$
by $v_1,\ldots,v_{n-1}$ where $w'(v_i) \leq w'(v_{i+1})$ for $i=1,\ldots,n-2$,
and where $w'$ is the sum of the labels given to edges incident with $v_i$.
Now assign the label $m-n+1+i$ to the edge $(v,v_i)$ for $i=1,\ldots,n-1$.
Notice that in the final labeling we have $w(v_i)=w'(v_i)+m-n+1+i$ and hence all
the weights of $v_1,\ldots,v_{n-1}$ are distinct. Also $w(v)=(n-1)(m-n+1)+n(n-1)/2$
and hence $w(v)$ is larger than any $w(v_i)$. Thus, $G$ is antimagic. \endpf

An {\em $S$-partial labeling} of $G=(V,E)$ is an assignment of distinct labels to
some of the edges of $G$ from a set of positive integers $S$ (we allow $|S| > |E|$ in the definition
of partial labeling). Given a partial labeling, the weight of a vertex $v$, denoted $w(v)$,
is the sum of the labels of the edges incident with $v$, that received a label.
A {\em completion} of an $S$-partial labeling is an assignment of distinct  unused labels from
$S$ to all the remaining non-labeled edges of $G$. We need the following simple lemma:

\begin{lemma}
\label{l42}
Let $G$ be any graph with $r$ vertices and $m$ edges.
Let $S$ be a set of positive integers with $|S|=m+2$. Then, any $S$-partial labeling that satisfies the property
that no more than $\lceil r/2 \rceil$ vertices have the same positive
weight, has a completion that also satisfies this property.
\end{lemma}
{\bf Proof}\,
Assume that $t \leq m$ edges are labeled. We need to label the remaining $m-t$
edges. We use induction on $m-t$. If $m-t=0$ there is nothing to prove.
Assuming the lemma holds when $m-t=k-1$, we prove it for $k$.
Assume, therefore, that $m-t=k > 0$.
Pick an arbitrary non-labeled edge, $e=(x,y)$. Since we have $m+2$ possible labels
and at most $m-1$ have been used, there are at least three possible labels
that we may assign to $e$. Denote the 
three non-used labels by $a_1,a_2,a_3$.
We will show that at least one of the three possible assignments maintains
the property in the statement of the lemma. If we assign $a_i$ to $e$ then,
after the assignment, we have that the total weight of $x$ is $w(x)+a_i$
and the total weight of $y$ is $w(y)+a_i$. The other vertices did not change their
total weight. Assume, for the sake of contradiction, that each of the three assignments
fails. This means that for each $i=1,2,3$, one of $w(x)+a_i$ or $w(y)+a_i$
appears at least $\lceil r/2 \rceil+1$ times as a total weight.
Hence, in the original partial labeling, we have that for each $i=1,2,3$, one of
$w(x)+a_i$ or $w(y)+a_i$ appears at least $\lceil r/2 \rceil$ times as a total weight.
Assume, without loss of generality, that $w(x)+a_1$ and $w(x)+a_2$ appear
at least $\lceil r/2 \rceil$ times as a total weight in the original labeling.
The original labeling also has the weight $w(x)$ appearing in $x$.
Thus, the overall number of vertices of $G$ is at least $2\lceil r/2 \rceil +1 > r$,
a contradiction. We have proved that we can label $e$ and maintain the desired property.
Now, we remain with only $k-1$ unlabeled edges, and we can complete the labeling using
the induction hypotheses. \endpf

Recall that an {\em even graph} is a graph whose vertices all have even degree.
It is easy and well-known that any even graph can be decomposed into
edge-disjoint simple cycles.
It is also easy to see
that any graph has a subforest such that the deletion of the
edges of this subforest from the graph results in an even graph.
We use these facts in the following two lemmas.
\begin{lemma}
\label{l43}
Let $G=(V,E)$ have $n$ vertices and $\Delta(G)=n-2$. If $|E|=m \geq 2n-4$
then $G$ is antimagic.
\end{lemma}
{\bf Proof}\
The statement implies $n \geq 4$. We can assume $n \geq 5$ since for $n=4$ the only possible graph satisfying
the assumption is $C_4$, which is trivially antimagic.
Let $v_n$ denote a vertex of maximum degree $n-2$.
Let $v_{n-1}$ be its unique non-neighbor and let $v_1,\ldots,v_{n-2}$ be the other vertices.
Consider the induced subgraph $G^*$ on all vertices except $v_n$.
$G^*$ has $n-1$ vertices. Let $F$ be a set of edges that induce a subforest in $G^*$ such that
the deletion of $F$ from $G^*$ results in an even subgraph $G'$ of $G^*$. Clearly, $|F| \leq (n-1)-1=n-2$.
Assign the edges of $F$ the smallest even weights $2,4,\ldots,2|F|$ arbitrarily. Notice that this can
be done since $2|F| \leq 2n-4 \leq m$. Now, partition the edges of $G'$ into edge-disjoint
simple cycles $C_1,\ldots,C_p$. We label the edges of $C_i$ sequentially from $i=1$ until $i=p$.
We first use all the remaining even weights, and when we exhaust them, we turn to using odd
weights, starting from the smallest odd weights and continuing sequentially.
Now, either all cycles have all their edges labeled with the same parity, or else at most one
$C_i$ contains both odd and even weights, such that the odd weights form one path in $C_i$
and the even weights form the remaining path in $C_i$. In the latter case, exactly two vertices
of $G^*$ end up with a total odd weight, and we can assume neither of them is $v_{n-1}$ since $C_i$
has at least three vertices, while in the former case, all vertices of $G^*$
end up with total even weight. We are now left with the largest $n-2$ odd labels for the use
in the edges incident with $v_n$. Denote the remaining labels by $t,t+2,\ldots,t+2n-6$
where $t+2n-6$ is either $m$ or $m-1$ (depending on the parity of $m$).
Assume, without loss of generality, that after labeling 
$G^*$, we have $w_{G^*}(v_1) \leq w_{G^*}(v_2) \leq \cdots \leq 
w_{G^*}(v_{n-2})$.
If all the $w_{G^*}(v_i)$ are even then after labeling $(v_i,v_n)$ with the label $t+2i-2$ for $i=1,\ldots,n-2$
we have that all the total weights of $v_1,\ldots,v_{n-2}$ are distinct and odd, while the total weight of $v_{n-1}$ stays even.
The total weight of $v_n$ is clearly the largest of all total weights.
Hence, the labeling is antimagic.
We may therefore assume that for precisely two indices $j$ and $k$, with $1 \leq j < k \leq n-2$
have $w_{G^*}(v_j)$ and $w_{G^*}(v_k)$ odd. Since $n-2 \geq 3$, there are at least three
odd labels $t,t+2,t+4$. We may assign two of them to $(v_n,v_j)$ and $(v_n,v_k)$ so as to guarantee
that the total weights of $v_j$ and $v_k$ after the assignment is even, but distinct from that of $v_{n-1}$ which
is also even, and distinct from each other. The third odd label not used from $t,t+2,t+4$, together with the other $n-5$
odd labels, are assigned, sequentially, to the other vertices from $\{v_1,\ldots,v_{n-2}\} \setminus \{v_j,v_k\}$.
Hence, the total weights of these vertices is odd and distinct, and the total weight of $v_n$ is the largest.
Hence, $G$ is antimagic.
\endpf

\begin{lemma}
\label{l44}
Let $G=(V,E)$ have $n \geq 4$ vertices and $\Delta(G)=n-2$. If $|E|=m \leq 2n-5$
then $G$ is antimagic.
\end{lemma}
{\bf Proof}\,
Again, we will assume $n \geq 5$ as the case $n=4$ is trivial to check.
As in the proof of 
Lemma \ref{l43}, let $v_n$ denote a vertex of degree $n-2$ and
let $v_{n-1}$ denote the unique non-neighbor of $v_n$.
If $v_{n-1}$ is an isolated vertex then we can use Lemma \ref{l41} for
the subgraph induced by all vertices except $v_{n-1}$ and obtain that $G$
is antimagic ($v_{n-1}$ will be the unique vertex with total weight $0$ in this case).
Thus, we assume $v_{n-1}$ is not isolated.
As in Lemma \ref{l43}, let $G^*$ be the subgraph induced by all vertices except $v_n$.
Let $s$ denote the number of edges of $G^*$. 
Hence, $s=m-(n-2) \leq (2n-5)-(n-2)=n-3$.
Notice that $s < m/2$.
Thus, we can assign all edges of $G^*$ only even weights.
We consider three cases.

If $m=2n-5$ then $s=n-3$ and we use
all the even weights for labeling $G^*$. Then, we use the $n-2$ odd weights
to label the edges incident with $v_n$. As in Lemma \ref{l43}, we assign the
labels so as to guarantee that all vertices $v_1,\ldots,v_{n-2}$ have distinct
total odd weight. $v_{n-1}$ has total even weight, and $v_n$ has maximum weight.
Hence, $G$ is antimagic.

If $m=2n-6$ or $m=2n-7$ we have only $n-3$ odd weights and $m-n+3=s+1$ even weights.
Assign arbitrarily the first $s-1$ even weights to all but one edge of $G^*$,
denoted $e=(x,y)$. We may assume $v_1 \notin \{x,y\}$ since $n \geq 5$ so $n-2 \geq 3$.
Let $r_1$ and $r_2$ denote the largest two even weights. We may choose one of them
for the label of $e$. Let $a_1$ denote the total weight of $v_1$ at this point and
let $a_2$ denote the total weight of $v_{n-1}$ at this point.
If $v_{n-1} \in \{x,y\}$ we select $r_1$ for the label of $e$ if
and only if $r_2+a_1 \neq a_2+r_1$. Otherwise we select $r_2$ for the label
of $e$ and notice that in this case we must have $r_1+a_1 \neq a_2+r_2$.
If $v_{n-1} \notin \{x,y\}$ we select $r_1$ for the label of $e$ if
and only if $r_2+a_1 \neq a_2$. Otherwise we  select $r_2$ for the label of $e$
and notice that in this case we must have $r_1+a_1 \neq a_2$.
In any case we have shown that we can select a label for $e$ such that if we select
the other label for $(v_n,v_1)$ we have that the total weight of $v_1$
is even and distinct from the total weight of $v_{n-1}$ which is also even.
The remaining $n-3$ edges incident with $v_n$ receive the odd labels so as to
guarantee that the final weights of $v_2,\ldots,v_{n-2}$ are odd and distinct.
Finally, $v_n$ has the maximum total weight, so $G$ is antimagic.

If $m \leq 2n-8$ we have at least $s+2$ even weights for the labeling.
Let $S$ denote the set of $s+2$ largest even labels.
Label one of the edges of $G^*$ incident with $v_{n-1}$ (recall that
$v_{n-1}$ is non-isolated) with the largest even label (this is either
$m$ or $m-1$, depending on the parity of $m$). 
By Lemma \ref{l42}
we can complete this labeling to a full labeling of $G^*$ which uses
only elements of $S$, such that no set of $\lceil (n-1)/2 \rceil +1$ vertices
in $G^*$ have the same positive total weight.
Now, consider the remaining $n-2$ weights. We know that $x$
of them are even and $n-2-x$ of them
are odd, and, trivially, $s+x \leq n-2-x \leq s+x+1$.
Thus, $x \leq (n-3)/2$ in any case.
Assume, without loss of generality, that $w_{G^*}(v_1) \leq w_{G^*}(v_2) \leq \cdots \leq w_{G^*}(v_{n-1})$.
Denote the remaining $x$ even labels by $\{r_1,\ldots,r_x\}$.
The odd labels are $1,3,5,\ldots,2n-5-2x$.
We assign the edge $(v_i,v_n)$ the even label $r_i$, for $i=1,\ldots,x$.
We assign the edge $(v_i,v_n)$ the odd label $2i-2x-1$ for $i=x+1,\ldots,n-2$.
There are two cases. If for all $i=1,\ldots,x$ we have $w_{G^*}(v_i)+r_i \neq w_{G^*}(v_{n-1})$
then, the final weights of $v_1,\ldots,v_x,v_{n-1}$ are even and distinct.
The final weights of $v_{x+1},\ldots,v_{n-2}$ are odd and distinct.
The final weight of $v_n$ is the largest of all, so $G$ is antimagic.
Hence, we may assume that for some $i$ we have $w_{G^*}(v_i)+r_i = w_{G^*}(v_{n-1})$.
Assume $i$ is minimal with this property.
We claim that $w_{G^*}(v_i) > 0$. Indeed, otherwise we would have $w_{G^*}(v_{n-1})=r_i$
but this is impossible since we have labeled one of the edges incident with $v_{n-1}$ with the largest
even label, which is larger than $r_i$.
Let $Z=\{j ~: ~ j \geq i ~,~ w_{G^*}(v_i)=w_{G^*}(v_j) \}$.
Recall that no set of $\lceil (n-1)/2 \rceil +1$ vertices
in $G^*$ have the same positive total weight. Thus, $|Z| \leq \lceil (n-1)/2 \rceil$.
Let $k$ denote the largest element in $Z$.
Notice that $k+x-i+1=|Z|+x \leq n-2$.
We modify the labeling of the edges $(v_j,v_n)$ for $j=i,\ldots,k+(x-i+1)$.
We label $(v_{k+j},v_n)$ with the even label $r_{i+j-1}$ for $j=1,\ldots,x-i+1$.
We label $(v_j,v_n)$ with the odd label $2(j-i)+1$ for $j=i,\ldots,k$.
Notice that now all the vertices $v_1,\ldots,v_{i-1},v_{k+1},\ldots,v_{k+x-i+1}$
receive distinct total even weights which are 
also distinct from the total weight of
$v_{n-1}$ which is the even number $w_{G^*}(v_{n-1})$. All the vertices $v_i,\ldots,v_k$
and $v_{k+x-i+2}, \ldots,v_{n-2}$
receive distinct total odd weights. 
Finally, $v_n$ receives the maximal weight. Hence,
$G$ is antimagic. \endpf

{\bf Proof of Theorem \ref{t2}:}\,
The proof follows immediately from Lemma \ref{l41},  Lemma \ref{l43} and Lemma \ref{l44}. \endpf

\section{Proof of Theorem \ref{t3}}
It is convenient to separate the proof of Theorem \ref{t3}
into two cases. The bipartite case and the $k$-partite case where $k \geq 3$,
as the proofs of these cases are different. It is also convenient to
prove the bipartite case using the equivalent matrix formulation.
The following two lemmas yield Theorem \ref{t3}.

\begin{lemma}
\label{l51}
For all $m+n > 1$, the cells of  an $m \times n$ matrix can be assigned the distinct integers
$1,\ldots,mn$ such that all $m+n$ rows and columns receive distinct sums.
\end{lemma}
{\bf Proof}\,
Assume, without loss of generality, 
that $m \leq n$. The case $m=1$ is trivial so we 
assume $2 \leq m \leq n$.
Furthermore, we may assume $n \geq 4$
since $K_{2,2}$, $K_{2,3}$ and $K_{3,3}$ are easily verified as antimagic.
We show how to assign the distinct numbers $1,\ldots,mn$ to the cells of a matrix $A_{m \times n}$.
We assign the numbers $(i-1)n+1,\ldots,in$ to the cells in row $i$, $i=1,\ldots,m$. The assignment
within each row is always monotone increasing if $i$ is odd or if $i$ is the last row.
Otherwise, it is monotone decreasing. Let $R(i)$ denote the sum of the elements in row $i$ and
let $C(j)$ denote the sum of the elements in column $j$. Clearly, $R(i)-R(i-1)=n^2$
for $i=2,\ldots,m$. Hence, all row sums are distinct and form an arithmetic progression with
difference $n^2$. Also, for all even $i < m$, the sum of the first $i$ rows is a constant vector.
Thus, if $m$ is odd we have $C(j)-C(j-1)=1$ for $j=2,\ldots,n$ and if $m$ is even then
$C(j)-C(j-1)=2$. In any case, $C(n)-C(1) \leq 2(n-1)$. Since $2(n-1) < n^2$ we have
that at most one column sum is equal to a row sum. Hence, we may assume
that for one specific pair, $C(j)=R(i)$. Clearly $i < m$ since the last row contains the $n$
largest elements and $n \geq m$. Assume first that $i > 1$. If $i$ is even then $A(i,1)-A(i-1,1)=2n-1$.
If $i$ is odd then $A(i,n)-A(i-1,n)=2n-1$. In any case we can replace the values of two adjacent
cells in rows $i$ and $i-1$ whose difference is precisely $2n-1$. Notice that the sums of
the columns do not change, and the sums of the rows, except $i$ and $i-1$, do not change.
The new sum of row $i$ is $R(i)-2n+1$ and the new sum of row $i-1$ is $R(i-1)+2n-1$.
However, $(R(i)-2n+1)-(R(i-1)+2n-1)=n^2-4n+2 > 0$ for $n \geq 4$. Hence, all row sums
are still distinct, but now $R(i)-2n+1$ is smaller than $C(1)$, so all row and column sums are distinct.
Assume next that $i=1$. If $m \geq 3$ we have $A(2,1)-A(1,1)=2n-1$
and when replacing them we have, as in the previous case, that all row sums remain distinct while the
new sum of the first row is $R(1)+2n-1$, which  is greater than $C(n)$. Finally, if $i=1$ and $m=2$
we can simply assign all odd numbers sequentially to the first row and all even numbers sequentially to the
second row. The largest column sum is $4n-1$ and the smallest row sum is $n^2$ and $n^2 > 4n-1$. \endpf

\begin{lemma}
\label{l52}
If $G$ is a $k$-partite graph with $k \geq 3$ then $G$ is antimagic.
\end{lemma}
{\bf Proof}\,
Let the sizes of the vertex classes be $n_1,\ldots,n_k$
where $n_i \leq n_{i+1}$ for $i=1,\ldots,k-1$.
We may assume $n_1 \geq 2$ since otherwise Lemma \ref{l41} applies.
Let $A$ denote a vertex class of size $n_1$
and let $B=V(G) \setminus A$ denote the remaining vertices.
Put $|B|=m$ and notice that $m \geq n_2+n_3 \geq 2n_2 \geq 2n_1 \geq 4$.
Let $q$ denote the number of edges with both endpoints in $B$.
Our initial labeling assigns arbitrary distinct labels from $1,\ldots,q$ to the edges with both
endpoints in $B$. We denote $B=\{u_1,\ldots,u_m\}$ where
$w'(u_i) \leq w'(u_{i+1})$ for $i=1,\ldots,m-1$ and $w'(u_i)$ is the weight of
$u_i$ after the initial labeling.
We now show how to complete the labeling by assigning the labels
$q+1,\ldots,q+mn_1$ to the edges incident with $A=\{v_1,\ldots,v_{n_1}\}$.
Let $c(v_i,u_j)$ denote the label received by $(v_i,u_j)$. We define:
$$
c(v_i ,u_j) =
\left\{
    \begin{array}{ll}
(i-1)m +j +q & j~odd
\\
(n_1 - i)m +j +q & j~even
\end{array}
\right.
$$
except for the case when $m$ is even, when we define $c(v_i , u_m ) = im+ q$.
The additional contribution to the weights of the vertices of $B$ is:
$$
{ n_1 \over 2 } \left(2q +2j +m(n_1 -1)\right).
$$
Since the last expression is an increasing function of $j$ we obtain that the final labeling
satisfies
$$
w(u_1) <w(u_2 ) < \cdots <w(u_m),
$$
as required.
Observe that,
$$
w(u_m ) \leq {n_1 \over 2} \left(2q +m(n_1 +1)\right) + 
q (m-n_2 ) - { { (m-n_2 )( m-n_2 -1)} \over {2}}.
$$
On the other hand it is easily verified that the values of $w(v_i)$
are:
$$
w(v_i ) = { {m} \over {2} } \left(2i + 2q + n_1 (m-1) \right)  ,~
for~ m~odd
$$
and
$$
w(v_i ) = { {m} \over {2} } \left( 4i+ 2q + n_1 (m-2) -1\right) 
,~
for~ m~even.
$$
In both cases we have $w(v_i ) < w(v_{i+1})$. Thus, it suffices to verify that $w(v_1) >w(u_m)$.
Following some simple calculations we obtain that we have to verify the inequality:
$$
{ m \over 2 } (m +mn_1 +2-n_1-2n_2) + { n_2 \over 2 } ( n_2 +1) +q(n_2 -n_1 ) > { mn_1 \over 2} (n_1 +2) .
$$
Indeed, this inequality is correct since already
$$
{ m \over 2 } (m +mn_1 +2-n_1-2n_2)  >{ mn_1 \over 2} (n_1 +2)
$$
as $m \geq 2n_2 \geq 2n_1 \geq 4$. \endpf

\section*{Acknowledgment}
The authors thank Yair Caro for the valuable correspondence concerning some important issues of the paper.


\begin{thebibliography}{99}

\bibitem{AlSp} 
N. Alon and J. H. Spencer, {\bf The Probabilistic Method},
{\em Second Edition},
Wiley, New York, 2000.


\bibitem{Bo} B. Bollob\'as, {\bf Extremal Graph Theory}, Academic Press,
London, 1978.


\bibitem{ErLo} P. Erd\"os and L. Lov\'asz,
{\em Problems and results on 3-chromatic hypergraphs and some related
questions}, Infinite and Finite Sets (A. Hajnal et al., eds.), North-Holland,
Amsterdam (1975), 609-628.

\bibitem{HaRi}
N. Hartsfield and G. Ringel,
Pearls in Graph Theory, 
{\it Academic Press, INC., Boston, 1990 (revised version, 1994),
pp. 108-109.}

\end{thebibliography}
\end {document}